\documentclass[10pt]{article}
\usepackage[utf8]{inputenc}
\usepackage[english,french]{babel}

\usepackage{amsmath}
\usepackage{amssymb}
\usepackage{amsthm}
\usepackage{bm}
\usepackage{fancyhdr}
\usepackage{subcaption}
\usepackage{graphicx}

\newcommand{\ud}{\mathrm{d}}

\newtheorem{theorem}{Theorem}[section]
\theoremstyle{definition}
\newtheorem{algorithm}[theorem]{Algorithm}

\DeclareMathOperator{\Lapl}{\triangle}

\fancypagestyle{firstpage}{
\fancyhf{}

}

\begin{document}
\title{A Discontinuous Galerkin like Coarse Space correction for
  Domain Decomposition Methods with continuous local spaces: the
  DCS-DGLC Algorithm}
\author{Kévin Santugini\thanks{Institut Mathématiques de Bordeaux}}
\maketitle
\selectlanguage{english}
\begin{abstract}
In this paper, we are interested in scalable Domain Decomposition
Methods (DDM). To this end, we introduce and study a new Coarse Space Correction
algorithm for Optimized Schwarz Methods(OSM): the DCS-DGLC algorithm.
The main idea is to use a Discontinuous Galerkin like formulation to
compute a discontinuous coarse space correction. 
While the local spaces remain continuous, the coarse space should be discontinuous
to compensate the discontinuities introduced by the OSM at the
interface between neighboring subdomains. The discontinuous coarse 
correction algorithm can not only be used with OSM but also be used
with any one-level DDM that produce discontinuous iterates. While ideas from Discontinuous
Galerkin(DG) are used in the computation of the coarse correction, the
final aim of the DCS-DGLC algorithm is to compute in parallel the discrete solution to the
classical non-DG finite element problem.
\end{abstract}

\selectlanguage{french}
\begin{abstract} 
Dans cet article, nous nous intéressons aux Méthods de Décomposition
de Domaines (DDM) scalables. À cette fin, nous introduisons et
étudions un nouvel algorithme,  le DCS-DGLC, de correction grossière pour les méthodes de
Schwarz optimisées.
L'idée principale est d'utiliser une formulation apparentée aux
méthodes de Galerkin discontinues pour calculer une correction
grossière discontinue. Alors même que les espaces locaux restent
continus, l'espace grossier est choisi discontinu afin de  
pouvoir compenser les discontinuités introduites par les OSM aux
interfaces entre sous-domaines voisins. Cet algorithme de correction grossière
discontinue peut être employé non seulement avec les OSM mais aussi
avec toute DDM de un niveau qui produit des itérées discontinues. Bien que l'algorithme s'inspire des méthodes de
Galerkin discontinues, le but final de l'agorithme est de calculer en
parallèle la solution discrete à la formulation éléments finis classique, sans DG.
\end{abstract}
\selectlanguage{english}

\section*{Introduction}
During the second half of the $20$th century, the continuous increase
in computing power was mainly due to increased sequential computing
power. Parallel computers were very expensive and were
rarely available to researchers. During the past decade, the situation
changed as chip makers found it difficult to keep 
increasing the available sequential computing power. To keep offering
continuously increasing computing power, they turned to parallelism by putting
multiple cores in a single CPU. Nowadays, almost every computer that
is sold contains a multicore CPU, even laptops, and is thus a parallel
machine. It is therefore very important to design algorithms that can
fully use this parallel computing power. At the same time, massively
parallel computers also became
increasingly affordable to the professional market. Instead of just a handful of nodes, massively parallel
computers can consist of hundreds, thousands or tens of thousands of
nodes. Parallel machines of up to ten thousands of nodes are now
offered off the shelf and no longer need to be custom built. While these
massively parallel computers remain, for the time being, still too expensive for the
hobbyist, they can now be afforded by moderately rich research institutions.
Because of the ever
increasing availability of massively parallel computing power,
it is no longer sufficient for algorithms to be parallel. To take
advantage of the enormous parallel computing power of these machines,  
parallel algorithms must also be scalable. Scalability is the property for an
algorithm to work well when more and more nodes are added. One usually
distinguishes two kind of scalability:
\begin{description}
\item[Strong Scalability] The property for an algorithm to
  finish in half the time for a given problem size if the number of computation nodes is doubled.
\item[Weak Scalability] The property for an algorithm to finish in
  the same amount of time if both the problem size and the number of
  computation nodes are doubled.
\end{description}
Strong Scalability is hard to achieve as when the number of nodes
tends to infinity, the number of operation per processor would need to
go to zero. In this paper, we are interested in the design of weakly scalable
Domain Decomposition Methods.

Domain Decomposition Methods (DDM) are a family of algorithms designed to parallelize the
computation of numerical solutions to Partial Differential Equations.
These methods are designed with non shared memory architectures in
mind. In Domain Decomposition Methods, the computation domain is
subdivided in subdomains. In this paper, we only consider 
domain decompositions with no overlap between subdomains.
In DDM, the interior equation is solved in parallel in each subdomain
with artificial boundary conditions. Artificial boundary conditions
are updated using information coming from the other subdomains. 
In one-level DDM, the update only uses information coming from the
neighboring subdomains. 
One-level DDM can work fine when only a handful of subdomains are
present. However, one-level DDM cannot be scalable. Since, the global
solution often depends on what happens on the global domain, there can be
no reasonable hope a one-level DDM will converge in less iterations
than the diameter of the connectivity graph of the subdomain
decomposition. This means that for $1D$ problems, or for 
$2D$ problems on a strip, we must iterate at least as many times as there are
subdomains. Scalable Domain Decomposition Methods should have a
convergence rate independent of the number of subdomains.

To get scalable Domain Decomposition Methods, some kind of global
information transfer is needed. The standard way of making a DDM
scalable is adding a coarse space to a pre-existing one-level DDM,
thus making it a two-level DDM.  The first use of
coarse spaces in Domain Decomposition Methods can be traced back to 
\cite{Nicolaides:1987:DeflationConjugateGradientApplicationBoundaryValueProblems}.
Coarse spaces are used to send information
globally: from any subdomain to all other subdomains. This global
information transfer can be used to ensure scalability of the new DDM.
Among well known DDM with coarse spaces are the two-level Additive Schwarz
method~\cite{Dryja:1987:AVS}, the FETI method~\cite{Mandel:1996:BDP},
and the balancing Neumann-Neumann methods~\cite{Mandel:1993:BalancingDomainDecomposition,Dryja:1995:SMN,Mandel:1996:CSM}.
See~\cite{Toselli:2004:DDM,Smith:1996:DPM} for complete analyses of such
methods. 
Coarse spaces are currently an active area of research, for example
for high contrast problems
\cite{Dolean.Nataf.Scheichl.Spillane:AnalysisTwo-LevelSchwarzMethods,Nataf:2011:CSC}.
Adding an effective coarse space to Optimized Schwarz
Methods (OSM), or to any one-level DDM that produce discontinuous
iterates, is also highly non trivial:
see~\cite{Dubois:2009:CBO},~\cite[chap.5]{Dubois:2007:OSM} for 
numerous numerical tests, and~\cite{Dubois:2012:TOS} for a
rigorous analysis of a special case. 

For domain decomposition methods that produce discontinuous iterates
at the interface between subdomains, it is advantageous 
that the coarse space be discontinuous at the interfaces between 
subdomains. In~\cite{Gander.Halpern.Santugini:2013:DD21-DCSDMNV}, Discontinuous Coarse Spaces (DCS)
were advocated. In the same paper, was also introduced the 
DCS-DMNV (DCS-Dirichlet Minimizer Neumann Variational) algorithm that computed
a discontinuous coarse correction for the Optimized  Schwarz Method
(OSM) in FEM-based discretizations. See~\cite{Gander:2006:OSM}  for a
description and analysis of OSM. For a similar
approach to adding an efficient coarse space to the Restricted Additive
Schwarz (RAS) algorithm, see also~\cite{Gander.Halpern.Santugini:2013:ANC}. See~\cite{Efstathiou:2003:WRA} for an
analysis of the one-level RAS algorithm. 

In this paper, we introduce and analyze another two-level Domain
Decomposition Method: the DCS-DGLC (DCS-Discontinuous Galerkin Like
Correction) algorithm. Like the DCS-DMNV algorithm introduced 
in~\cite{Gander.Halpern.Santugini:2013:DD21-DCSDMNV}, the DCS-DGLC
algorithm makes use of Discontinuous Coarse Spaces, is suitable to
Finite Element discretizations, and designed to not need Krylov acceleration
to converge. The coarse correction step iself is different. In particular, in the DCS-DGLC, the computation of the coarse correction 
is inspired by the ideas present in Discontinuous Galerkin and 
uses a penalization function in the variational formulation at the
coarse level. Even though ideas from Discontinuous Galerkin are used, the goal of
the DCS-DGLC algorithm is to converge to the discrete solution of
the non discontinuous Galerkin finite elements formulation. It is
probable the coarse space algorithm could be adapted to a true DG
formulation but this goes far beyond the scope of this paper. 
In this paper, we only consider the $(\eta -\Lapl)u=0$ equation, both OSM
and the DCS-DGLC algorithm can be adapted to more general elliptic operators.

In~\S\ref{sect:Rappels}, we remind the readers about some previously known
ideas, results or algorithms concerning Domain Decomposition Methods and coarse
spaces that are useful to understand the DCS-DGLC algorithm.
Then, in~\S\ref{sect:DCS-DGLC}, we state the DCS-DGLC algorithm and 
explain the motivations and the ideas behind this algorithm.
Then, in~\S\ref{sect:NumericalResults}, we present some numerical
results for both the iterative version
in~\S\ref{subsect:NumericalResultsIterative}
and the Krylov accelerated version in~\S\ref{subsect:NumericalResultsKrylov}.

\section{Optimized Schwarz Methods and Discontinuous Coarse spaces}\label{sect:Rappels}
In this section, we recall the formulation of Optimized Schwarz
Methods, % in \S\ref{subsect:OSM}, 
the ideas behind the use of
Discontinuous Coarse Spaces (DCS), % in~\S\ref{subsect:DCS},
and the DCS-DMNV \cite{Gander.Halpern.Santugini:2013:DD21-DCSDMNV} 
algorithm.% in~\S\ref{subsect:DMNV}. 

%\subsection{One-level Optimized Schwarz Methods}\label{subsect:OSM}
Let $\Omega$ be a bounded domain of $\mathbb{R}^d$. Let
$(\Omega_i)_{1\leq i\leq N}$ be a non overlapping domain decomposition
of $\Omega$. The one-level Optimized Schwarz Methods are defined by
\begin{algorithm}[One-level Optimized Schwarz]\hspace{1cm}\\
\begin{enumerate}
\item Set an initial $u_i^0$.
Until convergence
\begin{enumerate}
\item Set $u_i^{n+1}$ as the unique solution to
\begin{align*}
(\eta-\Lapl)u_i^{n+1}&=f\quad\text{in $\Omega_i$}\\
\mathcal{B}_{ij}u_i^{n+1}&=\mathcal{B}_{ij}u_j^{n}\quad\text{on $\partial\Omega_i\cap\partial\Omega_j$}\\
u_i^{n+1}&=0\quad\text{on $\partial\Omega_i\cap\partial\Omega$}\\
\end{align*}
\end{enumerate}
\end{enumerate}
\end{algorithm}

%\subsection{The importance of discontinuous coarse spaces}\label{subsect:DCS}
As explained in the introduction, one-level Optimized Schwarz Methods cannot be scalable.
To make them two-level and scalable, we need a coarse space. Such a coarse space should contain
discontinuous functions. The motivations behind the use of discontinuous coarse spaces can be found in 
details in~\cite{Gander.Halpern.Santugini:2013:DD21-DCSDMNV}. The basic idea
is that since many DDM, and in particular Opimized Schwarz Methods (OSM), introduce discontinuities at the interfaces
between subdomains, we need coarse functions with discontinuities also located at the interface between subdomains
to compensate those discontinuities. To understand why, consider a
generic iterative Coarse Space Correction algorithm:
\begin{algorithm}[Generic Coarse Space Correction
  Algorithm]\label{algo:GenericCSCorrection}\hspace{1cm}\\
\begin{enumerate}
\item Choose a coarse space $X$.
\item Initialize $u_i^{0}$, either by zero or using the coarse solution.
\item For $n\geq 0$ and until convergence
\begin{enumerate}
\item In each subdomain $\Omega_i$,
  compute the local iterates $u_i^{n+1/2}$ in parallel using the
  optimized Schwarz algorithm.
\item Compute a coarse correction $U^{n+1}$ belonging to
  a coarse space $X$.
\item Set the global iterates to
  $u_i^{n+1}:=u_i^{n+1/2}+U^{n+1}$.
\end{enumerate}
\item Set either $u_i:=u_i^{n-1/2}$ or  $u_i:=u_i^{n}$ where $n$ is the exit
  index of the above loop.
\end{enumerate}
\end{algorithm}
In the above algorithm, we haven't specified yet how to compute $U^{n+1}$.
More important than the algorithm used to compute $U^{n+1}$ is the
choice of the coarse space $X$ itself. A function $u$ in $\otimes_{i=1}^N H^1(\Omega_i)$ is
a weak solution to the linear elliptic equation $\eta-\Lapl u=f$ on $\Omega$ 
if and only if 
\begin{enumerate}
\item for all subdomains $\Omega_i$, $u_{\vert\Omega_i}$ is a weak solution to the interior equation $(\eta-\Lapl)u=f$ inside 
  subdomain $\Omega_i$.
\item there is no Dirichlet jumps, \textit{i.e.}, $u_i=u_j$ on $\partial\Omega_i\cap\partial\Omega_j$.
\item there is no Neuman jumps, \textit{i.e.}, $\frac{\partial
    u_i}{\partial\bm{n}_i}+\frac{\partial u_i}{\partial\bm{n}_i}=0$ 
on $\partial\Omega_i\cap\partial\Omega_j$.
\end{enumerate}
The coarse step in the algorithm should give global iterates $u_i^{n+1}$ that are
closer to satisfying these three conditions that the local iterates
$u_i^{n+1/2}$. Since the interior equation is already satisfied inside
each subdomain by the
local iterates, no improvement is possible there. The best that can be
achieved is for the global iterates to also satisfy the interior
equation inside each subdomains. Therefore, the coarse space elements
should satisfy the homogenous interior equation inside each
subdomains. If the coarse space $X$ was a subset of $H^1(\Omega)$, then the Dirichlet jumps of the global iterates across subdomains 
would be equal to the Dirichlet jumps of the local iterates and there
could be no improvement. If the global iterates are to have lower
Dirichlet jumps across subdomains than the local iterates, then the
coarse space $X$ must contain discontinuous functions. This is why, the
coarse space $X$ should always be a subset of 
\begin{equation*}
  \mathcal{A}=\{u\in H^{1,disc}_0(\Omega),\ \forall i,\  
    (\eta-\triangle) u_{\vert \Omega_i}=0\}.
\end{equation*}
The space $\mathcal{A}$ is small for one-dimensional problems but very
big for two or higher dimensional problems. hence, using the full
optimal theoretical coarse space $\mathcal{A}$ is unpractical in to or
higher dimensions. Using of a suspace $X$ of small
dimension as the coarse space is necessary.

%\subsection{The DCS-DMNV algorithm}\label{subsect:DMNV}
To our knowledge, the first discontinuous coarse space correction algorithm that
didn't need Krylov acceleration to converge for Optimized Schwarz Methods
is the DCS-DMNV (DCS- Dirichlet Minimizer Neumann Variational)
algorithm. See~\cite{Gander.Halpern.Santugini:2013:DD21-DCSDMNV} for 
a description and analysis. We reproduce the DCS-DMNV algorithm here:
\begin{algorithm}[DCS-DMNV]\label{algorithm:DCS-DMNV}
\begin{enumerate}
\item Choose a coarse space $X_d$. Set $X_c=X_d\cap H^1(\Omega)$.
\item Initialize $u_i^{0}$ by either zero or $u^0_{\vert\Omega_i}$ 
  where $u^0$ is the coarse solution.
\item Until convergence
\begin{enumerate}
\item Compute  in parallel the local iterates $u_i^{n+1/2}\in
  H^{1}(\Omega_{i})$ from the global iterates $u_i^n$
 using Optimized Schwarz.
% by
% \begin{subequations}
% \begin{align}
% \eta u_i^{n+1/2}-\triangle{u_i^{n+1/2}}&=f\quad\text{in $\Omega_i$},\\
% \frac{\partial u_i^{n+1/2}}{\partial n_i}
% +p u_i^{n+1/2}
% &=\frac{\partial u_{j}^{n}}{\partial n_i}+p u_{j}^{n}
% \quad\text{on $\partial\Omega_i\cap\Omega_j$},\\
% u_i^{n+1/2}&=0 \quad\text{on $\partial\Omega_i\cap\Omega$}.
% \end{align}
% \end{subequations}
\item Define a global $u^{n+1/2} \in  H^{1,disc}_0(\Omega)$ as $u_i^{n+1/2}$ in $\Omega_i$.   Set $U^{n+1}$ as the unique function in $X_d$ such that
\begin{subequations}\label{subeq:CoarseCorrectorDef}
\begin{equation}\label{eq:DefCoarseCorrectionDirichletJumpMinimizers}
q(u^{n+1/2}+U^{n+1})=\min_{v\in X_d}q(u^{n+1/2}+v),
\end{equation}
with $q(u)=\sum_{ij}\int_{\partial\Omega_i\cap\partial\Omega_j}\lvert u_i-u_j\rvert^2\ud\sigma(\bm{x})$, and satisfying
\begin{equation*}%\label{eq:DefCoarseCorrectionWeak}
\eta\int_\Omega U^{n+1}(x) v(x)\ud x
+\int_\Omega  \nabla U^{n+1} (x)\nabla v(x)\ud x
-
\sum_{i,j}\int_{\partial\Omega_i\cap\partial\Omega_j}
\left(\frac{\partial u_i^{n+1/2}}{\partial n_{i}}
+\frac{\partial u_j^{n+1/2}}{\partial n_{j}}
\right)v \ud\sigma,
\end{equation*}
for all test functions $v$ in $X_c$.
\end{subequations}
\item Set $u_i^{n+1}:=u_i^{n+1/2}+U^{n+1}$.
\end{enumerate}
\item Set $u:=u_i^{n-1/2}$ on $\Omega_i$ for each $i$ in $\{1,\ldots,N\}$.
\end{enumerate}
\end{algorithm}
In the DCS-DMNV algorithm
the discontinuous coarse correction was computed by choosing the unique coarse corrector that minimized 
the $L^2$ norm of the Dirichlet jump and satisfied the weak
formulation for coarse test functions. Test functions
needs to be  $H^1$ in the variational formulation. Besides, for
symmetric problems, it was advantageous to choose the intersection
between the coarse space and $H^1$ as the set of test
functions. This added the requirement that the coarse space contained
a sufficiently large ``continuous'' subset. In the next section, we
introduce a new algorithm,  the DCS-DGLC algorithm, in which the
coarse correction problem no longer requires test functions that
belongs to $H^1(\Omega)$.

\section{The DCS-DGLC algorithm}\label{sect:DCS-DGLC}
 Our goal is
to design another discontinuous coarse space correction algorithm for OSM that can
use discontinuous test functions.  Inspired by Discontinuous Galerkin
formulations, we remove the Dirichlet jump minimization. As in
Discontinuous Galerkin methods, a penalization parameter $q>0$ is
introduced in front of a boundary term that penalizes jumps across the interfaces 
between neighboring subdomains. 
All functions in the coarse space can now be used as test functions for
the coarse problems instead of just those that happen to also be $H^1$.
\begin{algorithm}[DCS-DGLC]
\begin{enumerate}
\item Choose a coarse space $X$.
%\item Choose two positive real numbers $p$ and $q$.
\item Initialize $u_i^{0}$ by either zero or $u^0_{\vert\Omega_i}$ 
  where $u^0$ is the coarse solution.
\item Until convergence
\begin{enumerate}
%\begin{subequations}
\item Compute  in parallel  the local iterates $u_i^{n+1/2}\in H^{1}(\Omega_{i})$
  from the global iterates $u_i^n$ using Optimized Schwarz.
%% begin{align}
% \eta u_i^{n+1/2}-\triangle{u_i^{n+1/2}}&=f\quad\text{in $\Omega_i$},\\
% \frac{\partial u_i^{n+1/2}}{\partial n_i}
% +p u_i^{n+1/2}
% &=\left(\frac{\partial u_{j}^{n}}{\partial n_i}+p u_{j}^{n}+\frac{\partial u_{j}^{n}}{\partial n_i}+p u_{j}^{n}\right)
% \quad\text{on $\partial\Omega_i\cap\Omega_j$},\\
% u_i^{n+1/2}&=0 \quad\text{on $\partial\Omega_i\cap\Omega$}.
% \end{align}
%\end{subequations}
\item Define a global $u^{n+1/2} \in  H^{1,disc}_0(\Omega)$ as
  $u_i^{n+1/2}$ in $\Omega_i$.   
\item Set $U^{n+1}$ as the unique function in $X$ such that
\begin{equation}\label{eq:DGCoarseCorrectionVariational}
\begin{split}
&\phantom{=}
\sum_{i=1}^N\eta\int_\Omega U^{n+1}(x) V(x)\ud x
+\sum_{i=1}^N\int_\Omega  \nabla U^{n+1} (x)\nabla V(x)\ud x\\
\\&\phantom{=}
+q\sum_{ij}\int_{\Gamma_{ij}}(u_i^{n+1/2}+U_i^{n+1}-u_j^{n+1/2}-U_j^{n+1})\cdot(V_i-V_j)\ud\sigma(\hat{\bm{x}})
\\&=-
\sum_{ij}\int_{\partial\Omega_i\cap\Omega_j}\theta
\left(\frac{\partial u_i^{n+1/2}}{\partial n_{i}}+\frac{\partial u_j^{n+1/2}}{\partial n_{j}}\right)(V_i+V_j)
\ud\sigma(\hat{\bm{x}}),
\end{split}
\end{equation}
for all test functions $V$ in $X$.
\item Set $u_i^{n+1}:=u_i^{n+1/2}+U^{n+1}$.
\end{enumerate}
\item Set $u:=u_i^{n-1/2}$ on $\Omega_i$ for each $i$ in $\{1,\ldots,N\}$.
\end{enumerate}
\end{algorithm}

As the DCS-DMNV algorithm, the DCS-DGLC is suitable for Finite Element based
discretizations. 
However, contrary to the DCS-DMNV algorithm, 
test functions need not be $H^1(\Omega)$. This is advantageous for
symmetric problems as it lifts the requirement that 
the continuous subset of the discontinuous coarse space 
has to be big enough.

 To choose the value of
$\theta$, one could get inspiration from the Neumann-Neumann methods that also
involve jumps of the Neumann boundary conditions, we could choose to set
$\theta=1/4$ on non crosspoints and $\theta=1/K^2$ at crosspoints
where $K$ equals the number of subdomains that meats at that
particular crosspoints,
see~\cite[ch. 6]{Toselli:2004:DDM}. For ease of 
implementation, we chose
$\theta=1/\max(K)^2$ where $\max(K)$ is the maximum over all crosspoints
of the number of subdomains that meet at that given crosspoint.

\section{Numerical Results}\label{sect:NumericalResults} 
In this section, we show convergence curves for the DCS-DGLC algorithm
when solving $-\Lapl u=f$
with homogenous Dirichlet conditions in $\lbrack0,4\rbrack^2$. We
consider $Q_1$ finite elements on a cartesian mesh. We have $5\times5$
subdomains and $50\times50$ cells per subdomains. We iterate on the errors:
we start with $f=0$ but also with random Robin boundary conditions on
each subdomain. This has the advantage, for the iterative DCS-DGLC, of
having the machine precision plateau given not by the machine epsilon but by its underflow level. This is
not the case for the GMRES accelerated version. The error curves
always show the $\log_{10}$ of the $L^\infty$ norm of the error.

Unfortunately, for ``historical reasons'', our implementation can only use one
particular piecewise harmonic discontinuous coarse
space. This space is constructed by taking linear Dirichlet boundary
conditions on each edge of a subdomain, then solving the homogenous
equation. In that particular run of tests, we chose $\eta=0$, so this gives a discontinuous $Q_1$ coarse
space.  This coarse space happens to have a large $H^1$ subset. Future implementations won't
have this limitation.

When using Optimized Schwarz on finite elements,  it is important that the
Robin boundary conditions be lumped,
see~\cite{Gander.Hubert.Krell:2013:OptimizedSchwarzDDFV}. Without
lumping the Robin boundary part of the rigidity matrix, we would
observe slow modes. 

\subsection{Iterative DCS-DGLC
  algorithm}\label{subsect:NumericalResultsIterative}

\begin{figure}[th]
% Simulations at
% Travail/Programmes/MPI/DCS_Eta-Lapl/DCS-FEM-1D-2D-2013-07-01-Krylov-lump-dg-overlump/SimulDG/ResultsRandomDMNV5x5sb50x50c10p
% Travail/Programmes/MPI/DCS_Eta-Lapl/DCS-FEM-1D-2D-2013-07-01-Krylov-lump-dg-overlump/SimulDG/ResultsRandomDMNV5x5sb50x50c15p
% Travail/Programmes/MPI/DCS_Eta-Lapl/DCS-FEM-1D-2D-2013-07-01-Krylov-lump-dg-overlump/SimulDG/ResultsRandomDMNV5x5sb50x50c1p
% Travail/Programmes/MPI/DCS_Eta-Lapl/DCS-FEM-1D-2D-2013-07-01-Krylov-lump-dg-overlump/SimulDG/ResultsRandomDMNV5x5sb50x50c20p
% Travail/Programmes/MPI/DCS_Eta-Lapl/DCS-FEM-1D-2D-2013-07-01-Krylov-lump-dg-overlump/SimulDG/ResultsRandomDMNV5x5sb50x50c2p
% Travail/Programmes/MPI/DCS_Eta-Lapl/DCS-FEM-1D-2D-2013-07-01-Krylov-lump-dg-overlump/SimulDG/ResultsRandomDMNV5x5sb50x50c5p
% Travail/Programmes/MPI/DCS_Eta-Lapl/DCS-FEM-1D-2D-2013-07-01-Krylov-lump-dg-overlump/SimulDG/ResultsRandomNoCoarse5x5sb50x50c10p
% Travail/Programmes/MPI/DCS_Eta-Lapl/DCS-FEM-1D-2D-2013-07-01-Krylov-lump-dg-overlump/SimulDG/ResultsRandomNoCoarse5x5sb50x50c15p
% Travail/Programmes/MPI/DCS_Eta-Lapl/DCS-FEM-1D-2D-2013-07-01-Krylov-lump-dg-overlump/SimulDG/ResultsRandomNoCoarse5x5sb50x50c1p
% Travail/Programmes/MPI/DCS_Eta-Lapl/DCS-FEM-1D-2D-2013-07-01-Krylov-lump-dg-overlump/SimulDG/ResultsRandomNoCoarse5x5sb50x50c20p
% Travail/Programmes/MPI/DCS_Eta-Lapl/DCS-FEM-1D-2D-2013-07-01-Krylov-lump-dg-overlump/SimulDG/ResultsRandomNoCoarse5x5sb50x50c2p
% Travail/Programmes/MPI/DCS_Eta-Lapl/DCS-FEM-1D-2D-2013-07-01-Krylov-lump-dg-overlump/SimulDG/ResultsRandomNoCoarse5x5sb50x50c5p
\begin{minipage}{\textwidth}
\begin{minipage}{0.45\textwidth}
\subcaption{One-level}\label{subfig:ConvergenceRatesNoCoarse}
\includegraphics[width=0.9\textwidth]{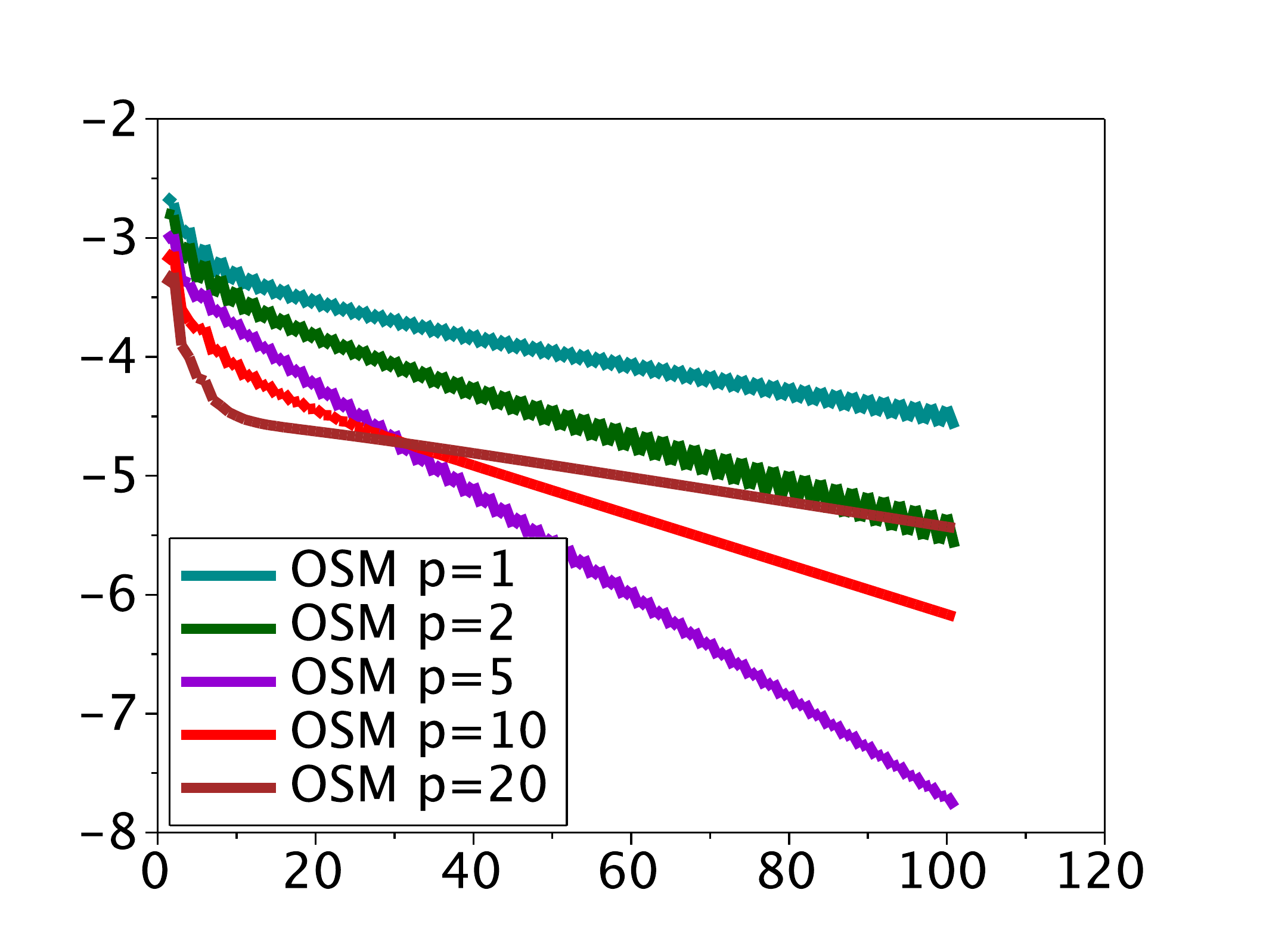}
\end{minipage}
\hspace{0.05\textwidth}
\begin{minipage}{0.45\textwidth}
\subcaption{DCS-DMNV}\label{subfig:ConvergenceRatesDMNV}
\includegraphics[width=0.9\textwidth]{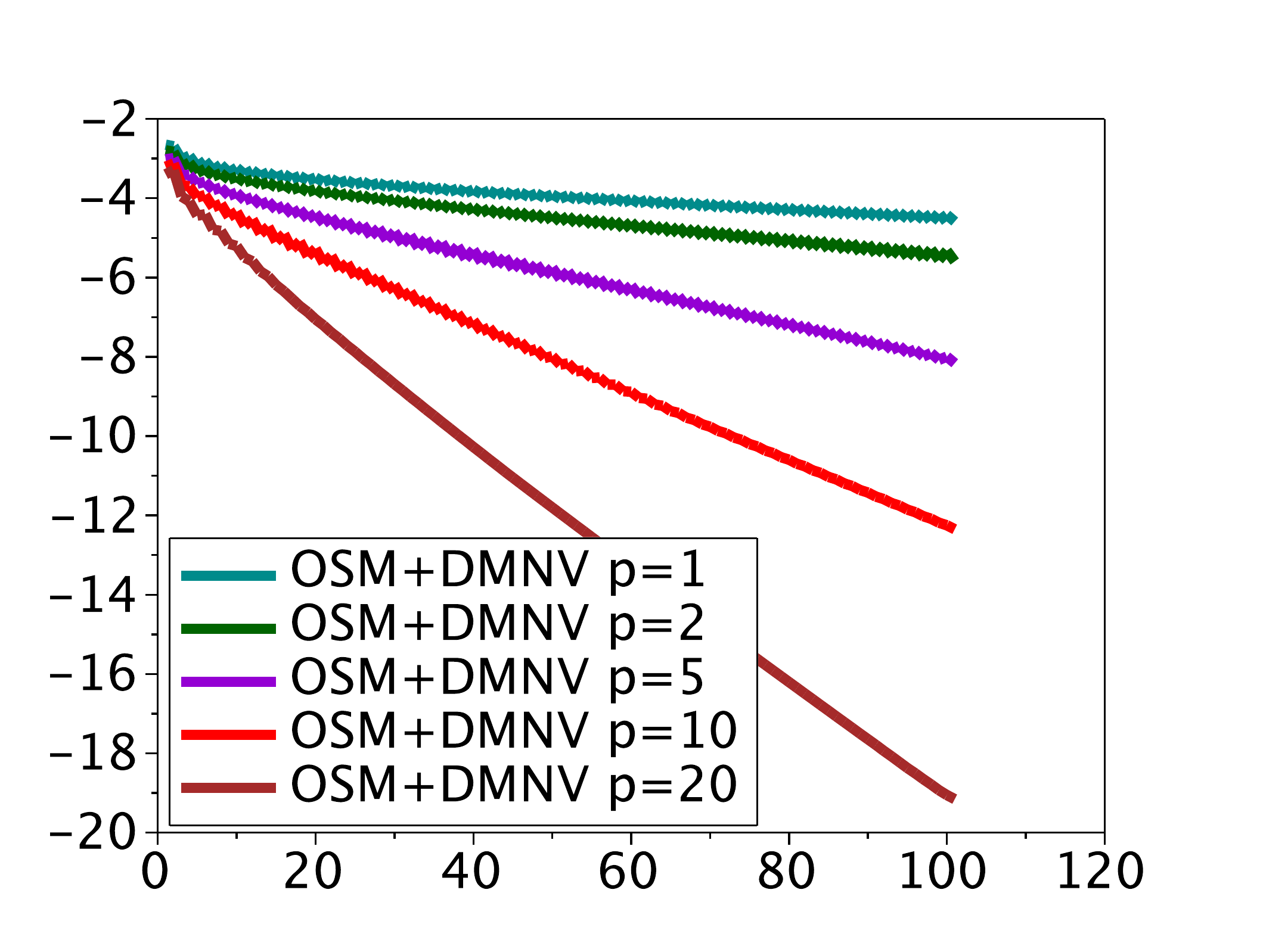}
\end{minipage}
\end{minipage}
\caption{Convergence rates for the one-level OSM and DCS-DMNV}\label{fig:ConvergenceRatesNoCoarseDMNV}
\end{figure}

\begin{figure}[th]
\begin{minipage}{\textwidth}
\begin{minipage}{0.45\textwidth}
\subcaption{$q=1$}\label{subfig:ConvergenceRatesDGLC1e0q}
\includegraphics[width=0.9\textwidth]{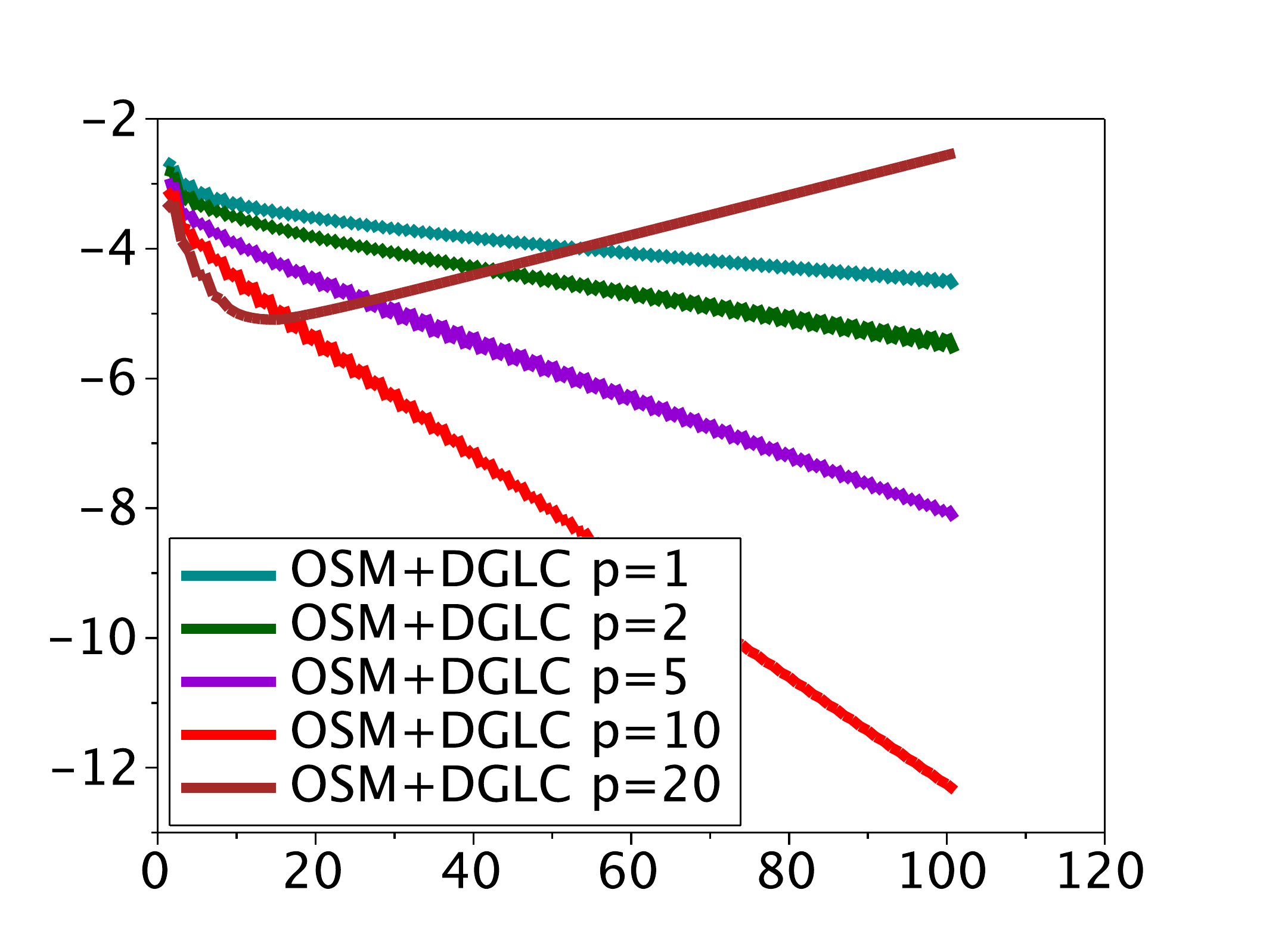}\\
\subcaption{$q=100$}
\includegraphics[width=0.9\textwidth]{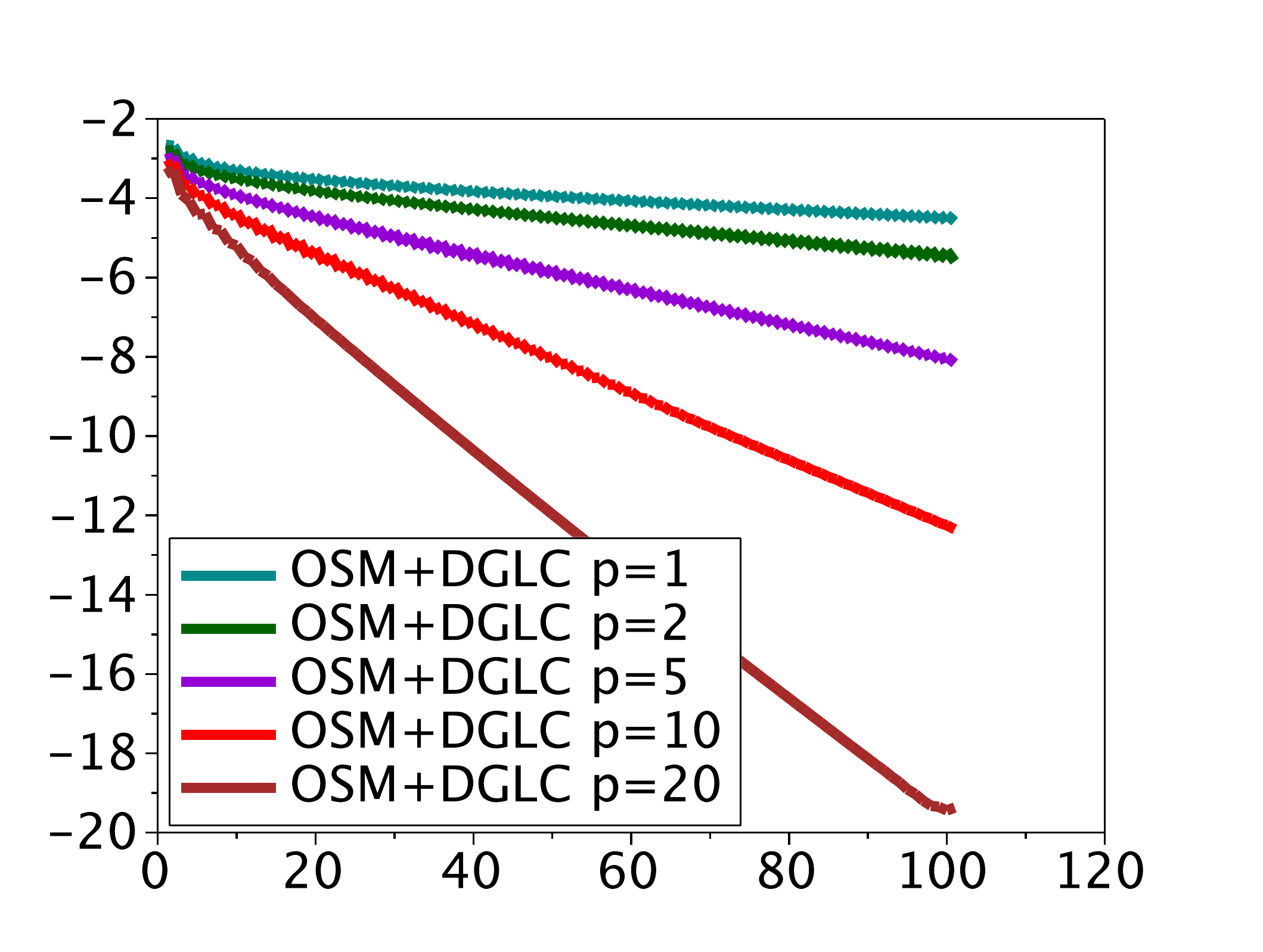}\\
\subcaption{$q=10000$}
\includegraphics[width=0.9\textwidth]{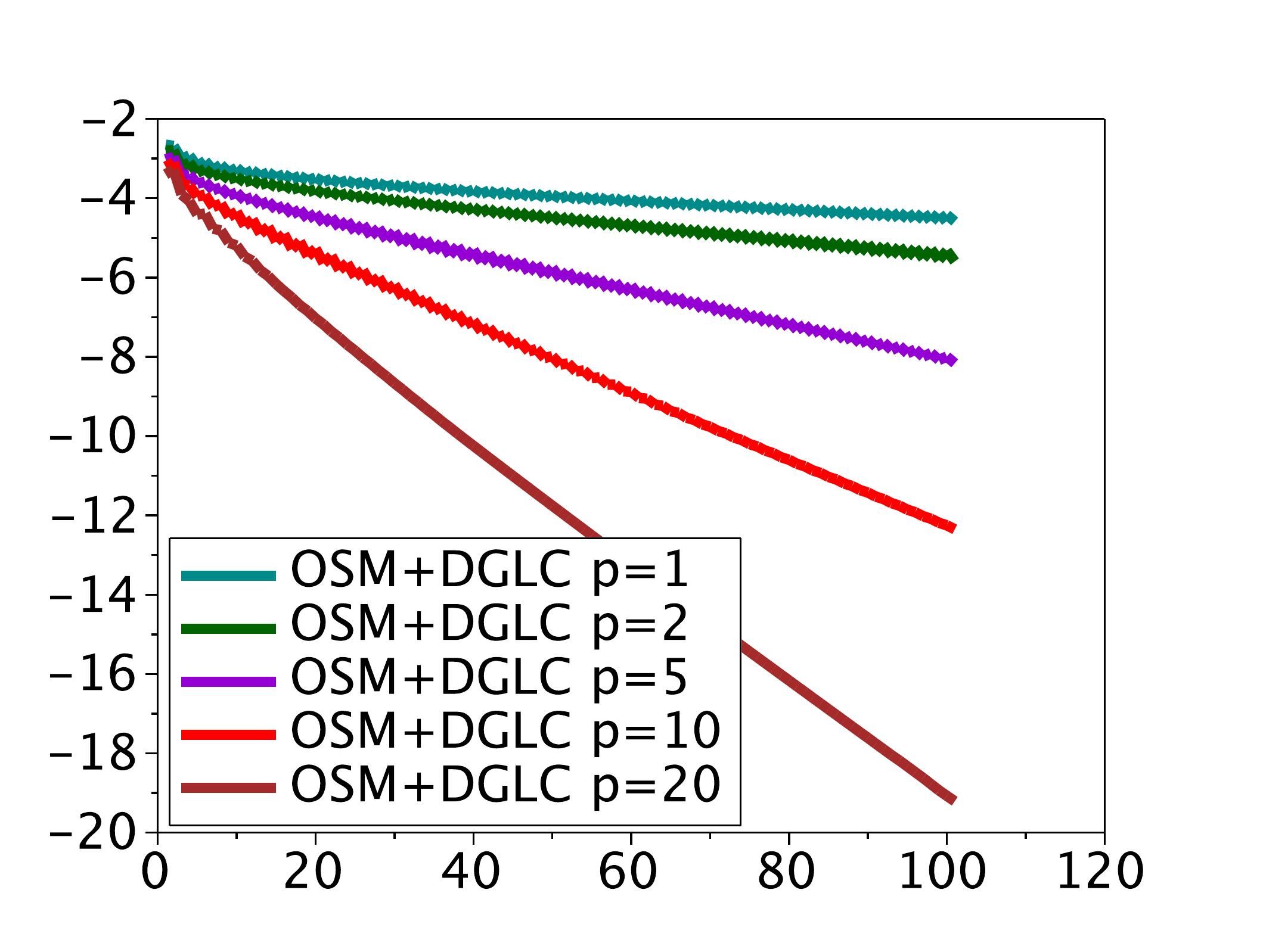}
\end{minipage}
\hspace{0.05\textwidth}
\begin{minipage}{0.45\textwidth}
\subcaption{$q=10$}
\includegraphics[width=0.9\textwidth]{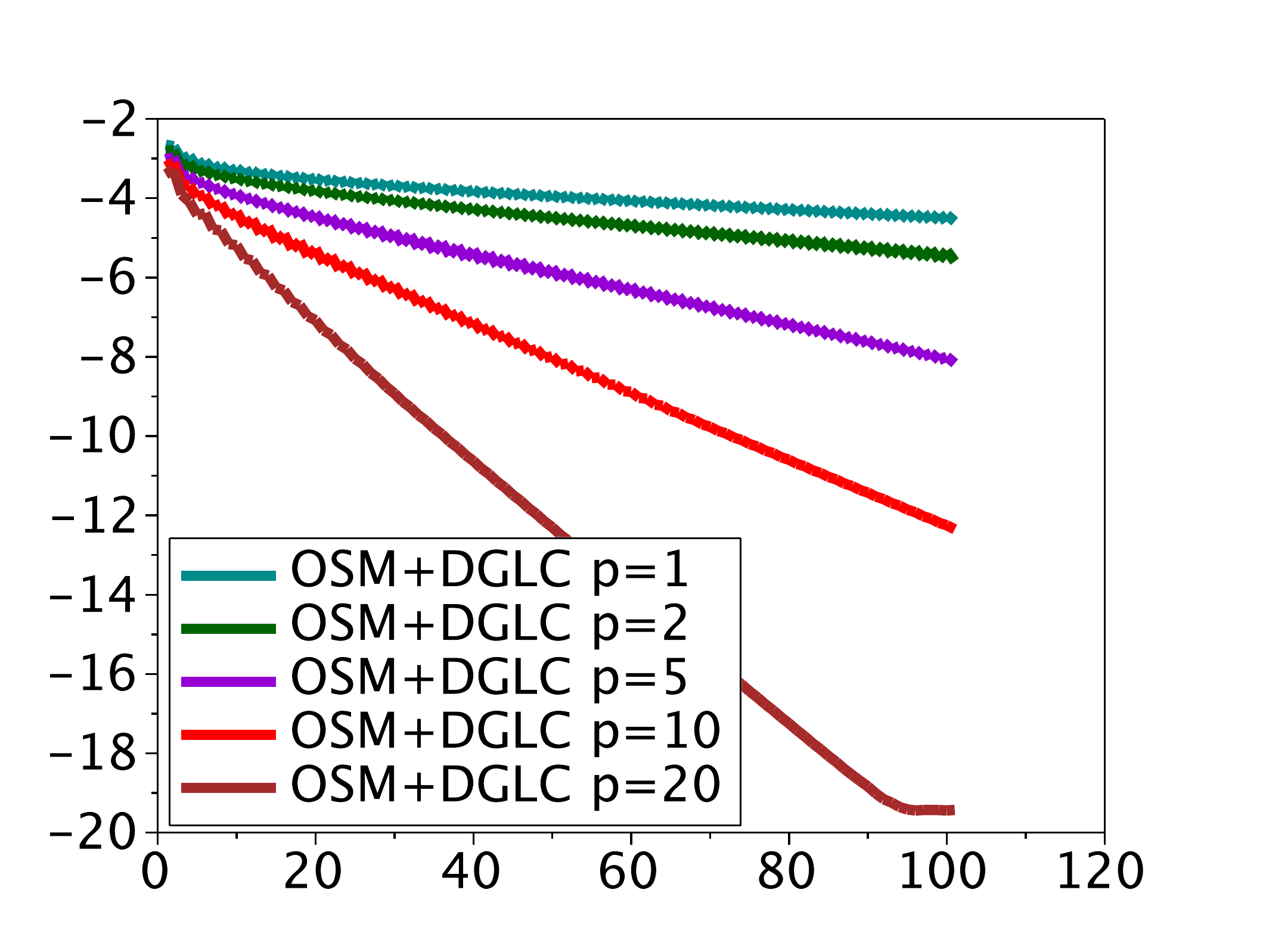}\\
\subcaption{$q=1000$}
\includegraphics[width=0.9\textwidth]{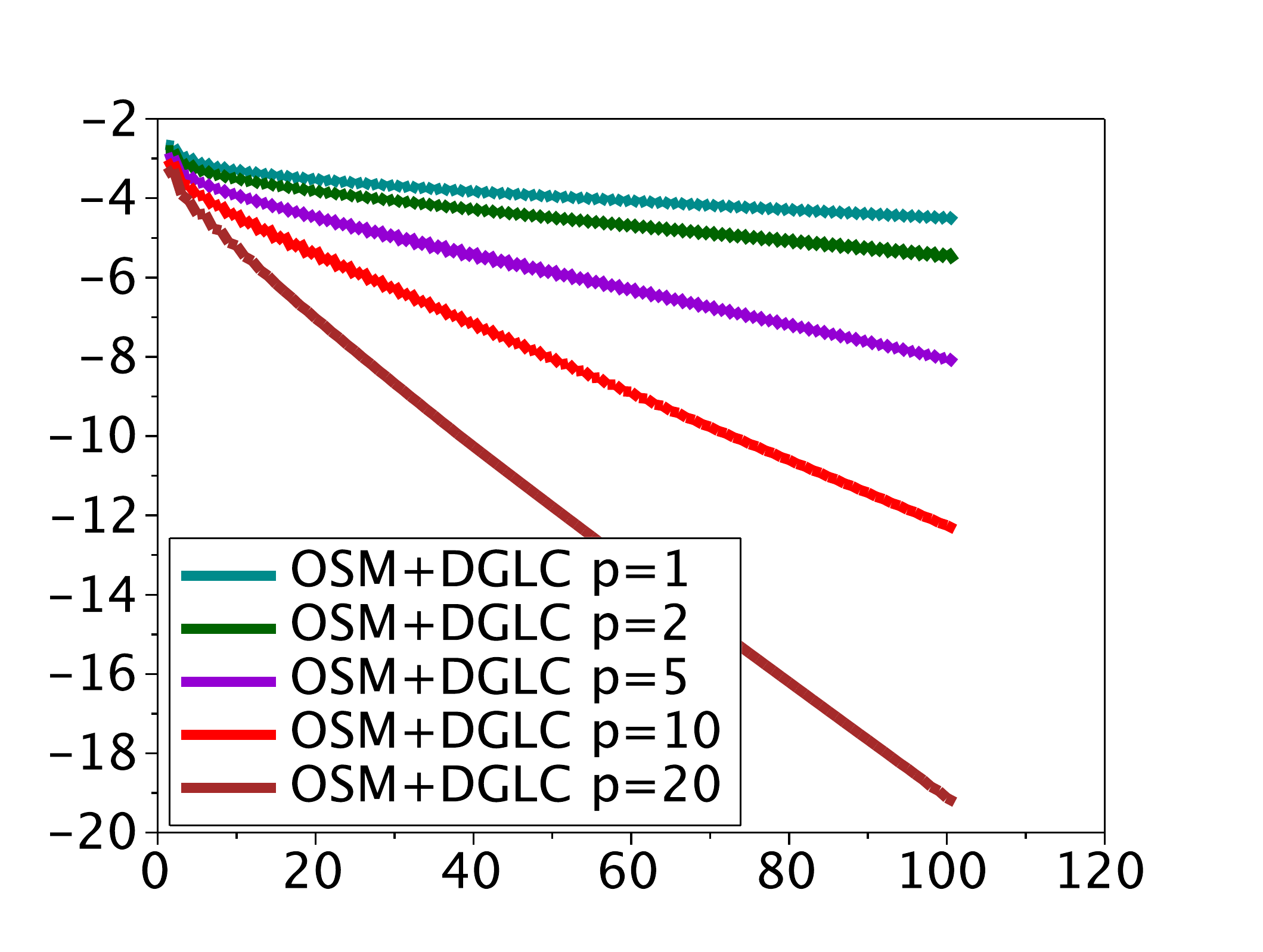}\\
\subcaption{$q=100000$}
\includegraphics[width=0.9\textwidth]{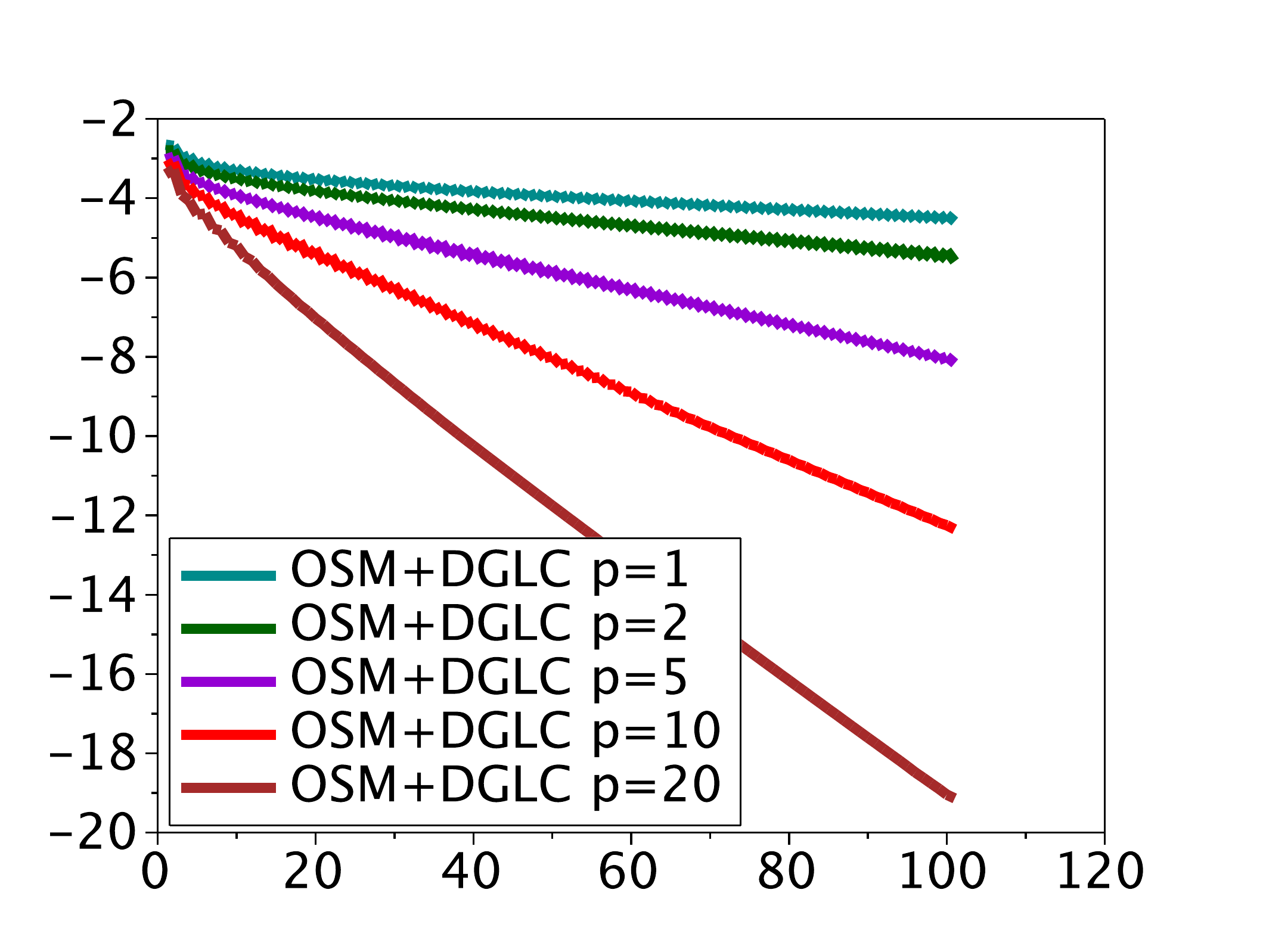}
\end{minipage}
\end{minipage}
\caption{Convergence rates for the DCS-DGLC Algorithm}\label{fig:ConvergenceRatesDGLC}
\end{figure}

In this section, we show error curves for the iterative DCS-DGLC
algorithm with six different values for the penalization parameter
$q$, see Figure~\ref{fig:ConvergenceRatesDGLC}. To give reference points for the performance of the
algorithm, we also show convergence curves for the one-level OSM and
the DCS-DMNV, see Figure~\ref{fig:ConvergenceRatesNoCoarseDMNV}. First
we observe that convergence is much slower for the one-level OSM, see
Figure~\ref{subfig:ConvergenceRatesNoCoarse}, which was to be
expected. We also observe that for $q=1$ and 
of $p=20$, the iterative DCS-DGLC  algorithm diverges, see Figure~\ref{subfig:ConvergenceRatesDGLC1e0q}.
For all the other values of $p$ and $q$, we observe convergence. The
performance of the DCS-DGLC is very close to the performance of the
DCS-DMNV algorithm, see Figure~\ref{subfig:ConvergenceRatesDMNV}. In
all cases both two-level algorithm converge much faster than the
one-level algorithm. For $p=20$, they reach an error of $10^{-20}$ in
$100$ iterations. We also observe that the behavior of the DCS-DGLC
algorithm seems to depend very little on $q$ once $q\geq10$.

\subsection{Krylov acceleration}\label{subsect:NumericalResultsKrylov}

\begin{figure}[ht]
\begin{minipage}{\textwidth}
\begin{minipage}{0.45\textwidth}
\subcaption{One-level}\label{subfig:ConvergenceRatesNoCoarse-GMRES}
\includegraphics[width=0.9\textwidth]{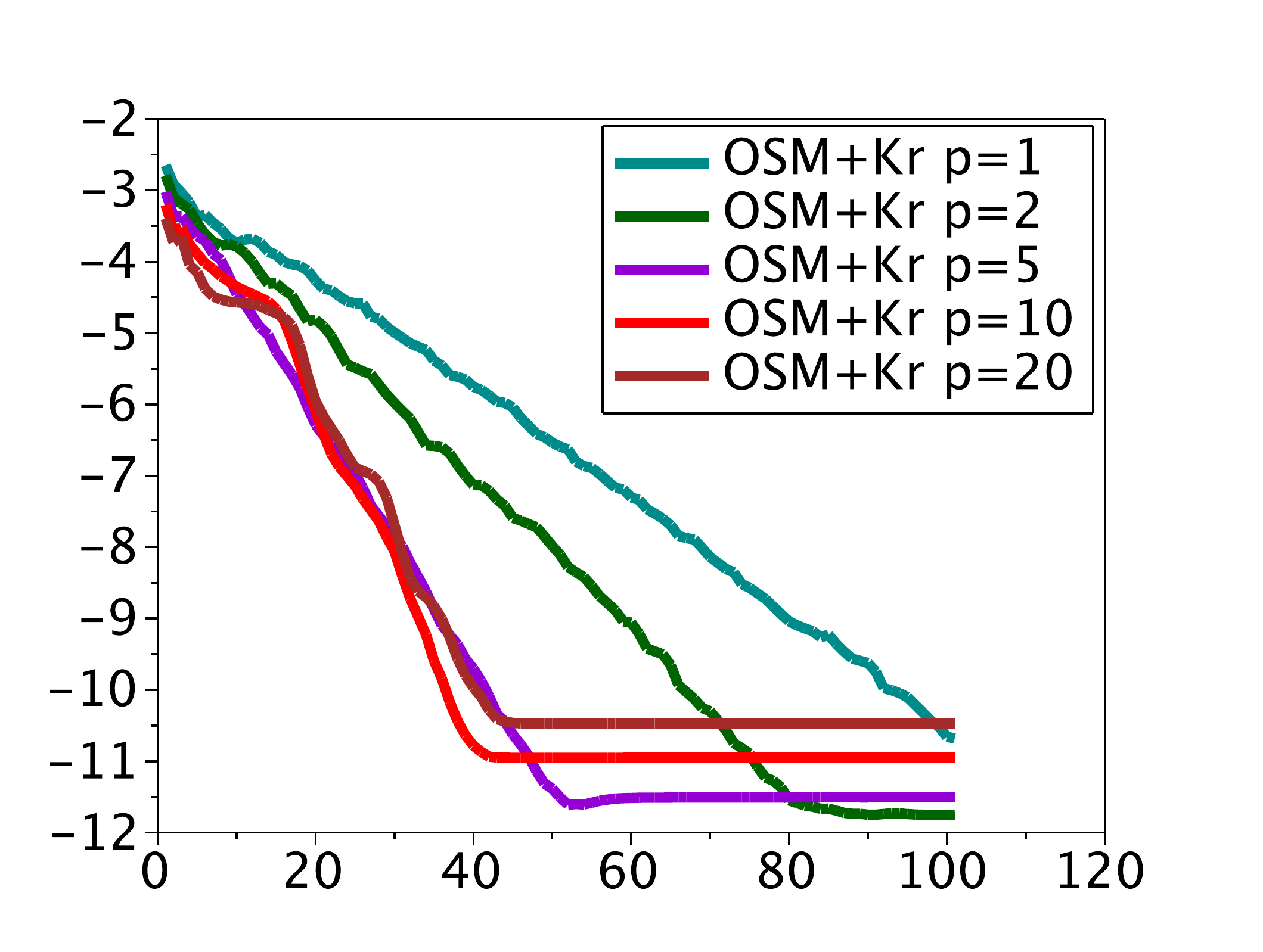}
\end{minipage}
\hspace{0.05\textwidth}
\begin{minipage}{0.45\textwidth}
\subcaption{DCS-DMNV}
\includegraphics[width=0.9\textwidth]{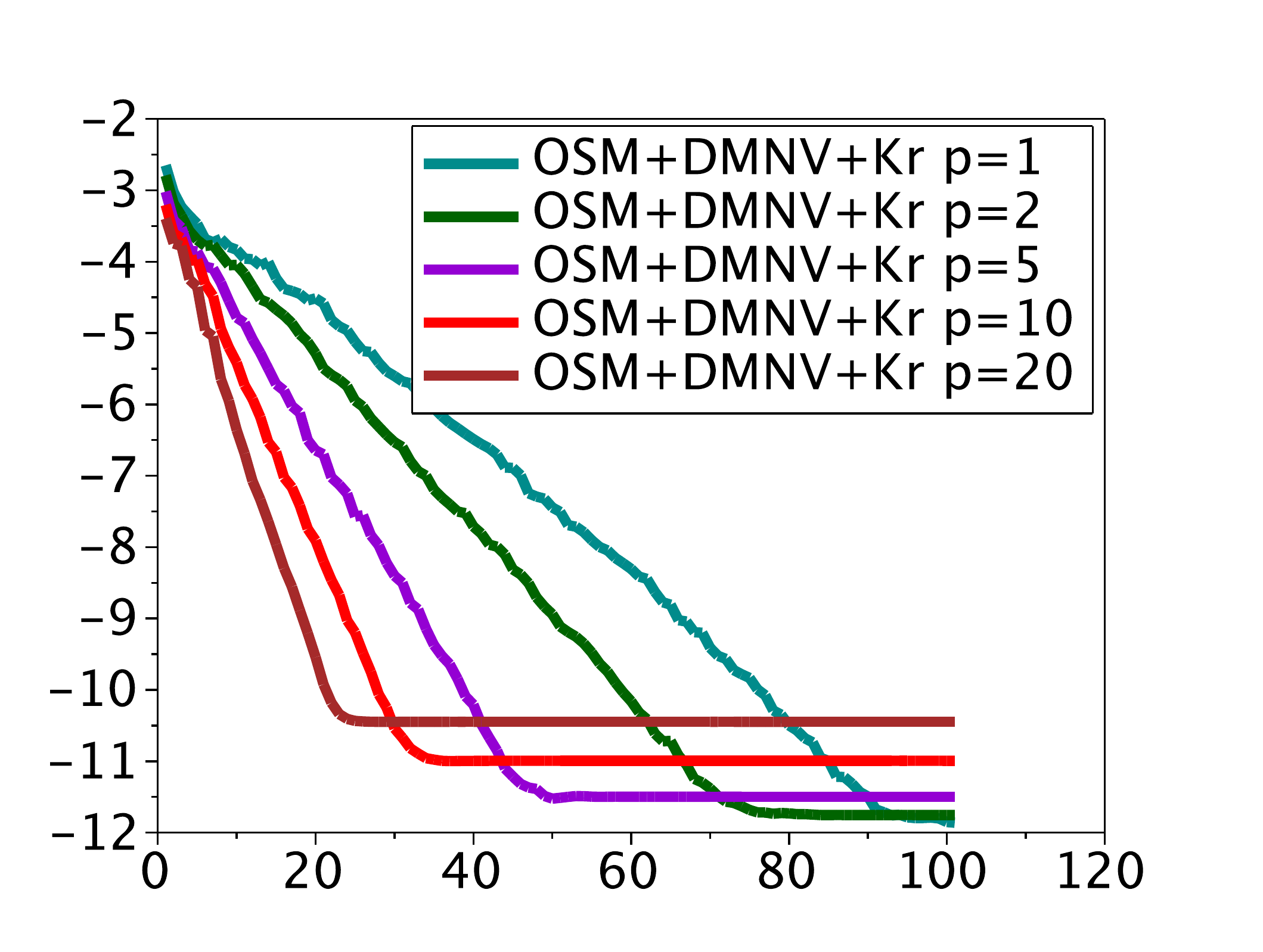}
\end{minipage}
\end{minipage}
\caption{Convergence rates for the one-level OSM and DCS-DMNV with
  GMRES acceleration}\label{fig:ConvergenceRatesNoCoarseDMNV-GMRES}
\end{figure}

\begin{figure}[th]
\begin{minipage}{\textwidth}
\begin{minipage}{0.45\textwidth}
\subcaption{$q=1$}
\includegraphics[width=0.9\textwidth]{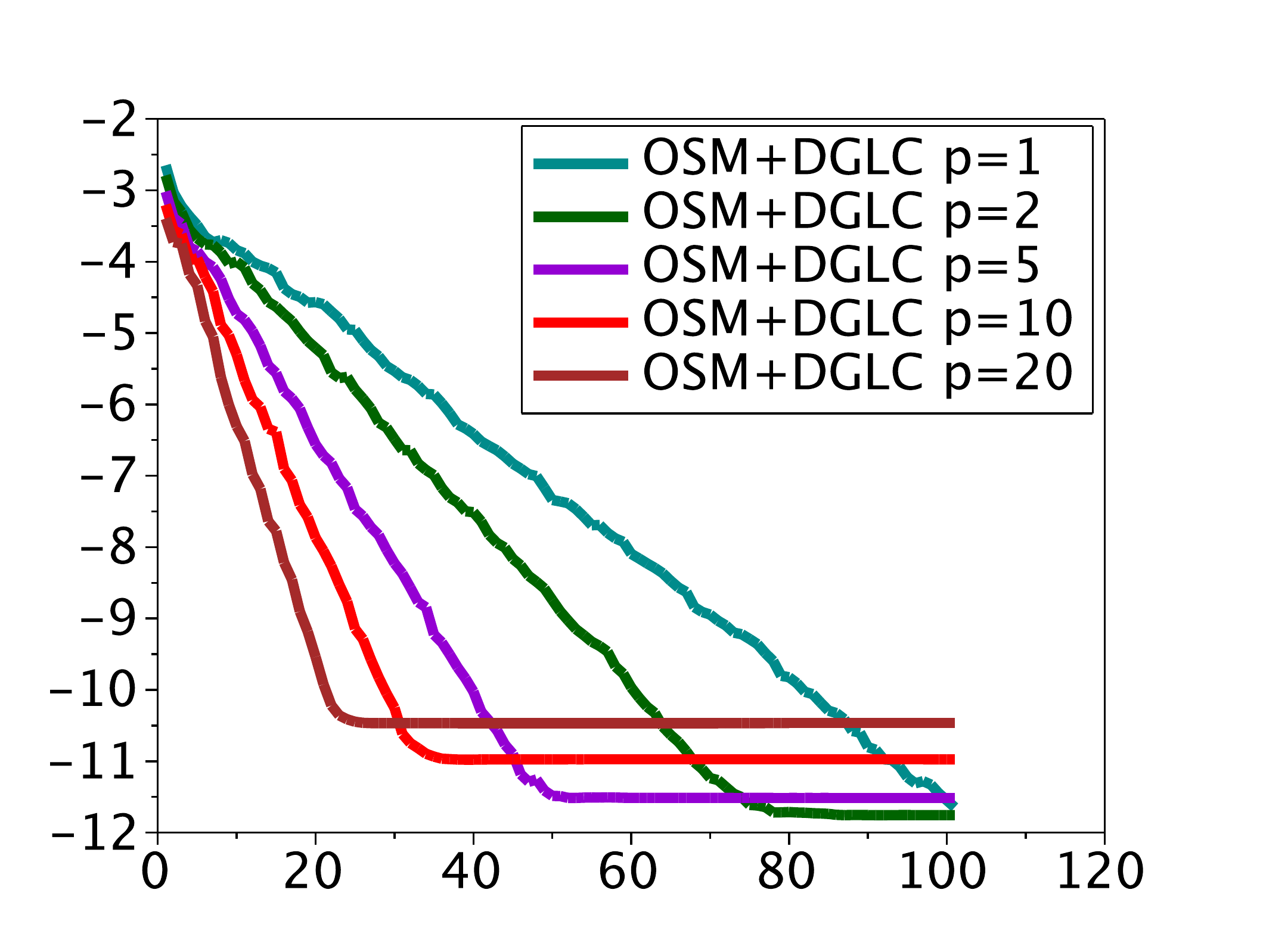}\\
\subcaption{$q=100$}
\includegraphics[width=0.9\textwidth]{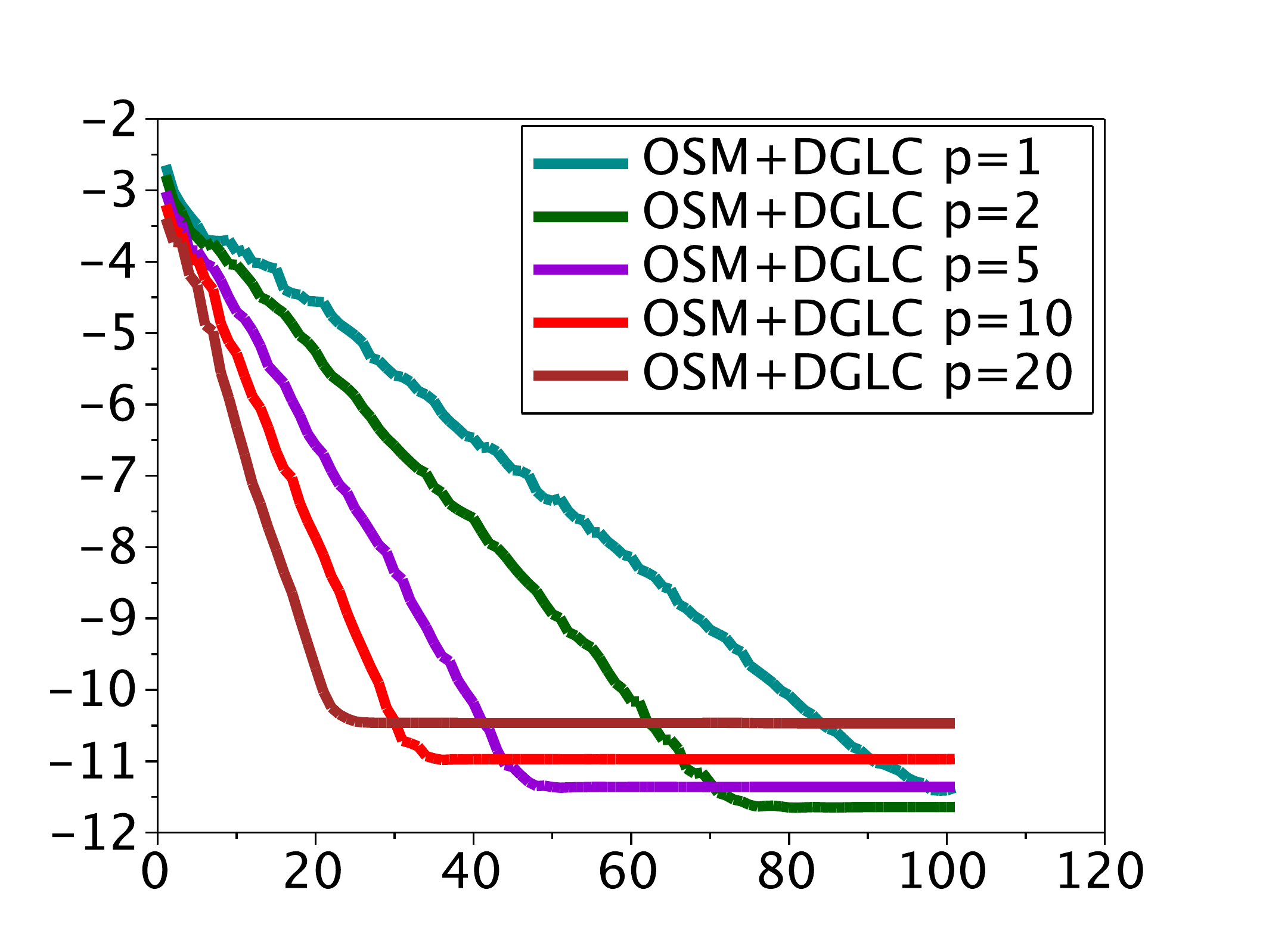}\\
\subcaption{$q=10000$}
\includegraphics[width=0.9\textwidth]{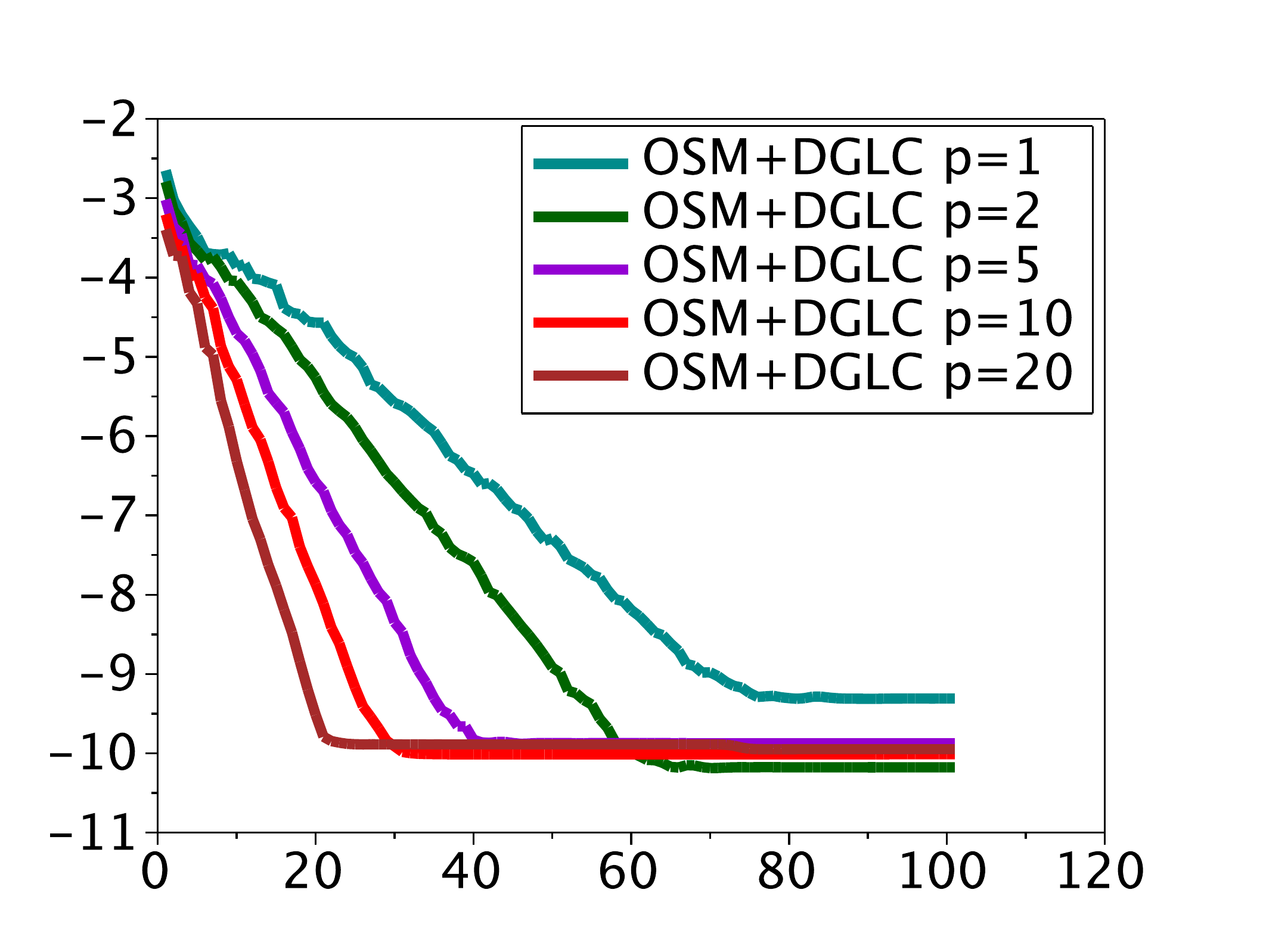}
\end{minipage}
\hspace{0.05\textwidth}
\begin{minipage}{0.45\textwidth}
\subcaption{$q=10$}
\includegraphics[width=0.9\textwidth]{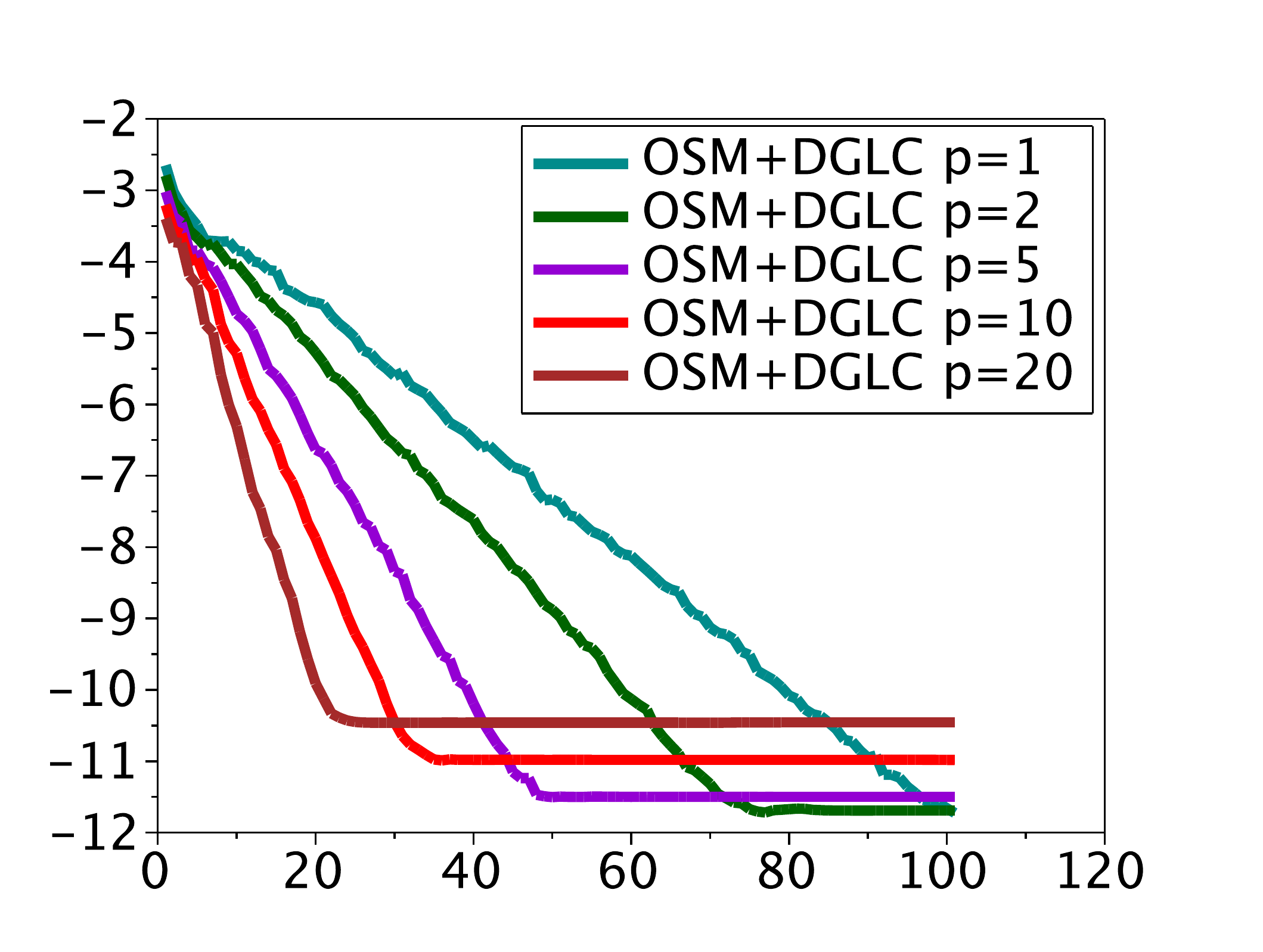}\\
\subcaption{$q=1000$}
\includegraphics[width=0.9\textwidth]{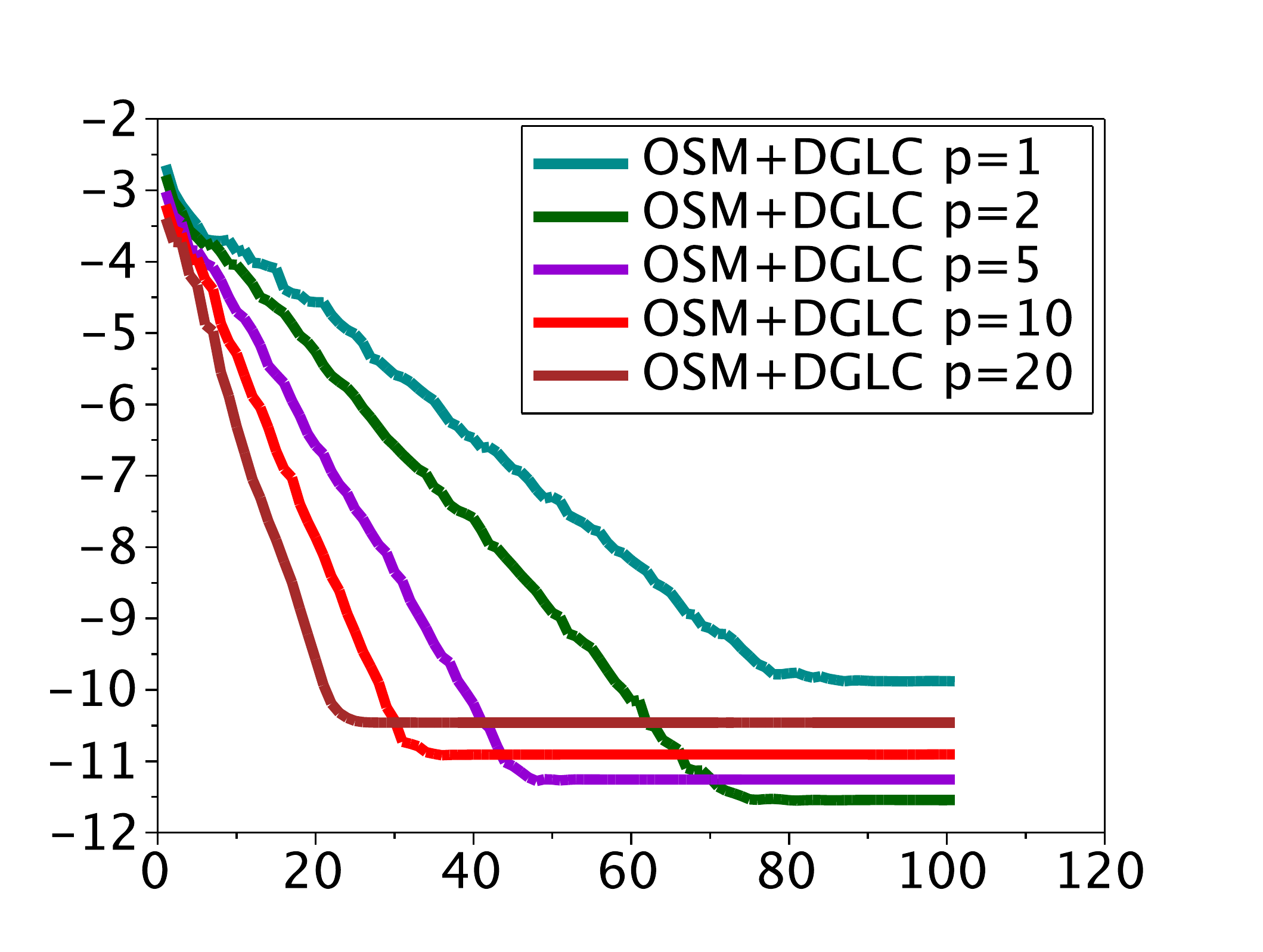}\\
\subcaption{$q=100000$}
\includegraphics[width=0.9\textwidth]{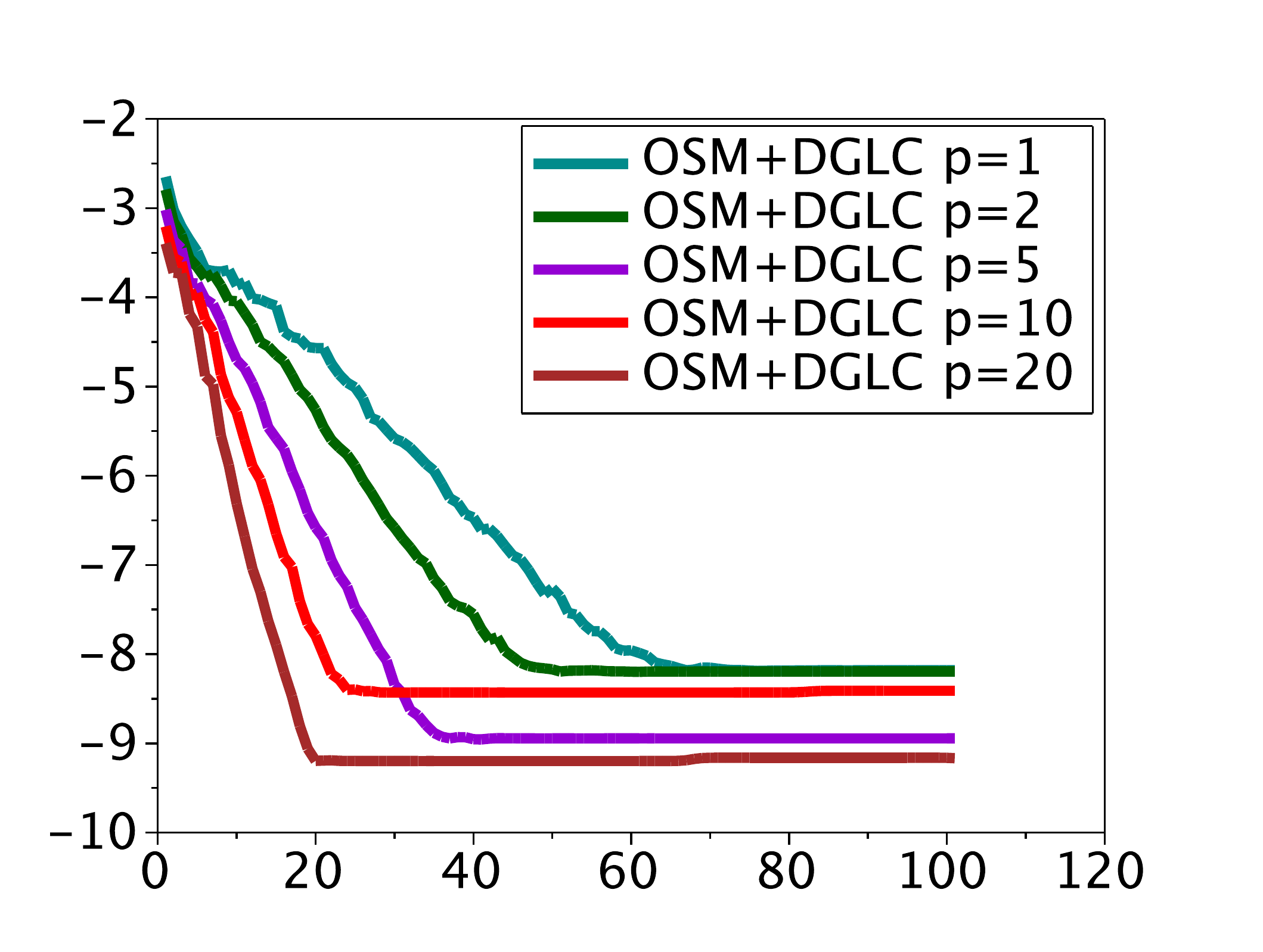}
\end{minipage}
\end{minipage}
\caption{Convergence rates for the DCS-DGLC Algorithm with GMRES acceleration}\label{fig:ConvergenceRatesDGLC-GMRES}
\end{figure}

It is well known that Domain Decomposition Methods can also be accelerated using 
Krylov methods. The classical idea that applies to any iterative method
is to see the iteration $u^{n+1}=Bu^{n}+b$ as a Richardson method to solve $(I-B)u=b$ and to
apply a Krylov method on $I-B$. While such acceleration gives faster methods, 
it is not sufficient to achieve scalable DDM. In practice, one
should always apply Krylov acceleration. However, for the purpose of
analysing an algorithm, it can be best to study the numerical behavior
of the iterative DDM itself. Krylov acceleration is so efficient it often hides away small 
design errors in a DDM algorithm. On the contrary, the slightest design
error often causes most iterative algorithms to fail.
This makes it easier to detect a DDM is non optimal and should be
improved. Nevertheless, we feel we wouldn't provide a complete picture
without providing results of numerical simulations in a Krylov setting.
It is important to note that Krylov acceleration cannot by itself make
a one-level DDM scalable.
% : in exact mathematics, Krylov methods are
% equivalent to extrapolation methods meaning the $n$th Krylov iterates
% belongs to TODO.

In this section, we show error curves for the GMRES accelerated DCS-DGLC
algorithm with six different values for the penalization parameter
$q$, see Figure~\ref{fig:ConvergenceRatesDGLC-GMRES}. To give reference points for the performance of the
algorithm, we also show convergence curves for the GMRES accelerated one-level OSM and
the GMRES accelerated DCS-DMNV, see
Figure~\ref{fig:ConvergenceRatesNoCoarseDMNV-GMRES}. Because of how GMRES
works, we reach a plateau that depends on the machine epsilon even
though we are iterating on the errors. The GMRES accelerated versions
converge much faster than their iterative counterparts. As expected,
the one-level method is the slowest. One advantage of the GMRES
accelerated DCS-DGLC is that it converges well even when $q=1$ and
$p=20$. The performance of the accelerated DCS-DGLC algorithm seems to
depend very little on $q$, however, the numerical plateau due to
rounding errors is higher for high values of $q$.

%  In the absence of Krylov acceleration, the slightest mistake in designing the algorithm can,
% and very often does, prevent convergence. This makes it easier to
% refine a DDM. Of course, in practice, once a satisfactory DDM
% has been designed and tested, it should always be accelerated 
% with Krylov methods.

\section*{Conclusion}
We have introduced a new discontinuous coarse space algorithm, the DCS-DGLC, that can
be used with any one-level Domain Decomposition Methods that produce
discontinuous iterates. We implemented that algorithm when used in
conjunction with Optimized Schwarz Methods, a subfamily of Domain
Decomposition Methods. Like its predecessor, the DCS-DMNV, which was
the subject of a previous paper, the DCS-DGLC is designed to work well
with Finite Element Discretizations. One potential advantage of the
DCS-DGLC over the DCS-DMNV is that it does not in theory needs that the 
coarse space contain a significant ``continuous''
subset. One advantage of lifting the 
``sizable continuous subset'' requirement is that it allows non
Dirichlet boundary conditions to be used to construct discontinuous coarse space
elements that also satisfy the interior equation inside each subdomain.
Unfortunately, due to the limitations of our implementation, 
we weren't able to study the numerical behavior of the
DCS-DGLC algorithm when the coarse space has a trivial ``continuous''
subspace. We plan to do so in the future.

\bibliography{ddm,ajoutmaths}

\begin{thebibliography}{10}

\bibitem{Dolean.Nataf.Scheichl.Spillane:AnalysisTwo-LevelSchwarzMethods}
Victoria Dolean, Fr\'ed\'eric Nataf, Robert Scheichl, and Nicole Spillane.
\newblock Analysis of a two-level schwarz method with coarse spaces based on
  local dirichlet to neumann maps.
\newblock {\em Computational Methods in Applied Mathematics}, 12(4):391--414,
  2012.

\bibitem{Dryja:1987:AVS}
Maksymilian Dryja and Olof~B. Widlund.
\newblock An additive variant of the {S}chwarz alternating method for the case
  of many subregions.
\newblock Technical Report 339, also Ultracomputer Note 131, Department of
  Computer Science, Courant Institute, 1987.

\bibitem{Dryja:1995:SMN}
Maksymilian Dryja and Olof~B. Widlund.
\newblock {S}chwarz methods of {N}eumann-{N}eumann type for three-dimensional
  elliptic finite element problems.
\newblock {\em Comm. Pure Appl. Math.}, 48(2):121--155, February 1995.

\bibitem{Dubois:2007:OSM}
Olivier Dubois.
\newblock {\em Optimized {S}chwarz Methods for the Advection-Diffusion Equation
  and for Problems with Discontinuous Coefficients}.
\newblock PhD thesis, McGill University, 2007.

\bibitem{Dubois:2009:CBO}
Olivier Dubois and Martin~J. Gander.
\newblock Convergence behavior of a two-level optimized {S}chwarz
  preconditioner.
\newblock In {\em Domain Decomposition Methods in Science and Engineering XXI}.
  Springer LNCSE, 2009.

\bibitem{Dubois:2012:TOS}
Olivier Dubois, Martin~J. Gander, Sebastien Loisel, Amik St-Cyr, and Daniel
  Szyld.
\newblock The optimized {S}chwarz method with a coarse grid correction.
\newblock {\em SIAM J. on Sci. Comp.}, 34(1):A421--A458, 2012.

\bibitem{Efstathiou:2003:WRA}
Evridiki Efstathiou and Martin~J. Gander.
\newblock Why {R}estricted {A}dditive {S}chwarz converges faster than
  {A}dditive {S}chwarz.
\newblock {\em BIT Numerical Mathematics}, 43(5):945--959, 2003.

\bibitem{Gander:2006:OSM}
Martin~J. Gander.
\newblock Optimized {S}chwarz methods.
\newblock {\em SIAM J. Numer. Anal.}, 44(2):699--731, 2006.

\bibitem{Gander.Halpern.Santugini:2013:ANC}
Martin~J. Gander, Laurence Halpern, and K{\'e}vin Santugini.
\newblock A new coarse grid correction for {RAS}.
\newblock In {\em Domain Decomposition Methods in Science and Engineering XXI}.
  Springer LNCSE, 2013.

\bibitem{Gander.Halpern.Santugini:2013:DD21-DCSDMNV}
Martin~J. Gander, Laurence Halpern, and Kévin Santugini-Repiquet.
\newblock Discontinuous coarse spaces for dd-methods with discontinuous
  iterates.
\newblock In {\em Domain Decomposition Methods in Science and Engineering XXI}.
  Springer LNCSE, 2013.

\bibitem{Gander.Hubert.Krell:2013:OptimizedSchwarzDDFV}
Martin~J. Gander, Florence Hubert, and Stella Krell.
\newblock Optimized {S}chwarz algorithm in the framework of {DDFV} schemes.
\newblock In {\em Domain Decomposition Methods in Science and Engineering XXI}.
  Springer LNCSE, 2013.
\newblock submitted.

\bibitem{Mandel:1993:BalancingDomainDecomposition}
Jan Mandel.
\newblock Balancing domain decomposition.
\newblock {\em Communications in Numerical Methods in Engineering},
  9(3):233--241, mar 1993.

\bibitem{Mandel:1996:BDP}
Jan Mandel and Marian Brezina.
\newblock Balancing domain decomposition for problems with large jumps in
  coefficients.
\newblock {\em Math.\ Comp.}, 65:1387--1401, 1996.

\bibitem{Mandel:1996:CSM}
Jan Mandel and Radek Tezaur.
\newblock Convergence of a {S}ubstructuring {M}ethod with {L}agrange
  {M}ultipliers.
\newblock {\em Numer.\ Math.}, 73:473--487, 1996.

\bibitem{Nataf:2011:CSC}
Fr\'ed\'eric Nataf, Hua Xiang, Victorita Dolean, and Nicole Spillane.
\newblock A coarse sparse construction based on local {D}irichlet-to-{N}eumann
  maps.
\newblock {\em SIAM J. Sci. Comput.}, 33(4):1623--1642, 2011.

\bibitem{Nicolaides:1987:DeflationConjugateGradientApplicationBoundaryValueProblems}
Roy~A. Nicolaides.
\newblock Deflation conjugate gradients with application to boundary value
  problems.
\newblock {\em SIAM J. Num. An,}, 24(2):355--365, 1987.

\bibitem{Smith:1996:DPM}
Barry~F. Smith, Petter~E. Bj{\o}rstad, and William Gropp.
\newblock {\em Domain Decomposition: Parallel Multilevel Methods for Elliptic
  Partial Differential Equations}.
\newblock Cambridge University Press, 1996.

\bibitem{Toselli:2004:DDM}
Andrea Toselli and Olof Widlund.
\newblock {\em Domain Decomposition Methods - Algorithms and Theory}, volume~34
  of {\em Springer Series in Computational Mathematics}.
\newblock Springer, 2004.

\end{thebibliography}
\bibliographystyle{plain}
\end{document}